\documentclass[12pt]{article}%
\usepackage{amsfonts}
\usepackage{amsmath}
\usepackage{amssymb}

\newcommand{\RR}{\mathbb{R}}
\newtheorem{theorem}{Theorem}[section]

\newtheorem{lemma}[theorem]{Lemma}

\newcommand{\qand}{\quad \mbox{and} \quad}

\newcommand{\pd}[2]{\frac{\partial #1}{\partial #2}}

\newcommand{\pro}{{\bf Proof. \hspace{2mm}}}

\begin{document}

\title{{\Large Necessary Optimality Conditions for
Some Control Problems of Elliptic Equations with Venttsel Boundary
Conditions}}
\author{Y. Luo \\
School of Mathematical and Geospatial Sciences, \\RMIT University,
GPO Box 2476V \\ Melbourne, Vic. 3001, AUSTRALIA \\
 email: yluo@rmit.edu.au}

\date{}
\maketitle

\begin{abstract}
In this paper we derive a necessary optimality condition for a
local optimal solution of some control problems. These optimal
control problems are governed by a semi-linear Vettsel boundary
value problem of a linear elliptic equation. The control is
applied to the state equation via the boundary and a functional of
the control together with the solution of the state equation under
such a control will be minimized. A constrain on the solution of
the state equation is also considered.

\end{abstract}



\section{Introduction}
\label{intro}
In this paper we discuss the necessary optimality conditions for a
class of optimal boundary control problems governed by a linear
elliptic partial differential equation with nonlinear Vettsel
boundary condition. We formulate the problem first.

Let $J :  C^{\alpha} (\partial \Omega) \rightarrow \RR$ be the
objective functional defined by
    \begin{equation}\label{cost} J(u) =\int_{\Omega} f(x, y_u(x)) \; dx
    +\int_{\partial \Omega} g(s, y_u(s),u(s)) \; ds
    \end{equation}
where $f : \Omega \times \RR \rightarrow \RR$ and $g :
\partial \Omega \times \RR \times \RR \rightarrow \RR$ are of class $C^1$ and
$y_u = G(u)$ is the solution of the state equation
    \begin{equation}\label{state1}  -\Delta y + y =0
     \quad
    \mbox{in}   \quad  \Omega, \qquad \qquad
     B(x, y, u)=0  \quad
    \mbox{on}   \quad   \partial \Omega.
    \end{equation}
corresponding to $u$. The boundary condition $B(x, y, u)=0 $ in
the state equation will be defined in details in Section 2.   Let
$U$ be a bounded set of functions, called the set allowable
controls, in $ C^{\alpha}(\partial \Omega)$ and $F :
C^{2,\alpha}(\overline{\Omega}) \rightarrow \RR$ be a constrain
functional on the state $y_u $ given by
     \begin{equation}\label{constrain} F(u) =\int_{\Omega} a(x, y_u(x)) \; dx
    +\int_{\partial \Omega} b(s, y_u(s)) \; ds
    \end{equation}
where $a : \Omega \times \RR \rightarrow \RR$ and $b :
\partial \Omega \times \RR  \rightarrow \RR$ are both of class
$C^1$. The function spaces $ C^{\alpha}(\partial \Omega)$ and $C^{2,\alpha}(\overline{\Omega}) $ will also be introduced in details in Section 2.  The control problem is formulated as follows:
    \begin{equation}\label{prob}  \qquad \left\{ \begin{array}{l}  \mbox{Minimize} \; \; J(u) \\
    u(x) \in U \\
    F(y_u) = 0
    \end{array}
    \right.
    \end{equation}
When the state equation is an elliptic equation with a traditional
boundary condition, i.e. either a Dirichlet or a Neumann boundary
condition, or a combination of the two in the form of a general
oblique boundary condition, this problem has been well studied.  A
survey of those results is given in \cite{ben}. For first and
second order necessary optimality conditions for state equations
with Neumann boundary conditions we refer to \cite{bon} and
\cite{casas}.

However problem (\ref{prob}) has not been studied so far when the
boundary condition is of Vettsel type.  A Vettsel boundary
condition consists of not only the unknown function and its first
order derivatives but also the second order tangential
derivatives of the unknown function. It has been shown in \cite{ven} that Venttsel
boundary condition is the most general admissible boundary
condition for second order elliptic operators.

Such a boundary condition can also be found in many engineering
problems and we refer \cite{apush:1} and \cite{luo:1} for details.
A simple example of it is the problem of heat conduction in a
medium enclosed by a thin skin and the conductivities of the
medium and the surrounding skin are significantly different, see
\cite{cj}. Generally speaking, all physical phenomena involving a
diffusion process along the boundary will give rise to a Vettsel
type boundary condition.

In the following, we will first discuss the existence theory of
the state equation in Section 2. In Section 3 We will establish
the differentiability of all functionals associated to our problem
and find the derivatives of them. In Section 4, we will state and
prove the main theorem of this paper. Finally in Section 5 we will
make some comments on work in progress and further development.

\section{State Equation} \label{sect:2}
Let $\Omega$ be a bounded open subset of $\RR^n$ with a $C^3$
boundary $\partial \Omega$. Let $\varphi : \partial \Omega \times
\RR \times \RR \rightarrow \RR$ be a $C^1$ function. Given a
function $u \in C^{\alpha}(\partial \Omega)$ we consider the
following boundary value problem:
\begin{eqnarray} \label{state} -\Delta y + y =0
     \qquad &
    \mbox{in}  & \quad  \Omega, \nonumber \\
     \Delta_{\partial \Omega} y +
    \partial_{\nu} y = \varphi (\cdot, y,u)  \quad &
    \mbox{on}  & \quad   \partial \Omega
    \end{eqnarray}
The boundary condition in (\ref{state}) is a special case of the
general Venttsel boundary condition.  The definition of such a
boundary condition is given as follows. Let
$\nu=(\nu^1,\ldots,\nu^n)$ be the inward unit normal vector field
on $\partial\Omega$.  Then the inward normal derivative of $u$,
denoted by $\partial_{\nu} y$, is defined by
\[ \partial_{\nu}y = D y \cdot \nu \]
where $Dy$ is the gradient vector.  Now we define the tangential
differential operators. Let $\{c^{ik}\}_{n\times n}$ be the matrix
whose entries are given by
\[c^{ik} =\delta^{ik}-\nu^i \nu^k\]
where $ \delta^{ik} $ is the Kronecker symbol. Then the first and
the second order tangential differential operators are then
defined by
\[\partial_i =c^{ik}D_k,\quad \partial_{ij} =\partial_{i} \partial_{j}
,\quad i,j,=1,\cdots,n,\] hence the tangential gradient operator
is defined by
\[\partial =(\partial_1, \cdots,\partial_n).\]
In particular the Laplace-Beltrami operator on the boundary
manifold is then defined by
\[ \Delta _{\partial \Omega} = \partial_{i} \partial_{i}. \]
Note that the second order tangential derivatives so defined are
not symmetric in general.

In order to understand the optimal control problem, the
fundamental issues is the existence and uniqueness of solutions of
the state problem (\ref{state}), as well as the continuous
dependence of the solutions upon the input $u$. The existence and
uniqueness of solutions of linear problems have been studied in
\cite{luo:1}. Similar results for quasi-linear equations with
quasi-linear boundary conditions are covered in \cite{luo:2}. Our
problem (\ref{state}) is a linear equation with a semi-linear
boundary condition, so the solvability and uniqueness can be
deduced from the general frame work of \cite{luo:2}. In order to
make the article more readable, without having to verify those
complicated general structure conditions stated in \cite{luo:2},
we prove the existence and uniqueness using only the results of
\cite{luo:1}.

We seek classical solutions in the H\"{o}lder space
$C^{2,\alpha} (\overline{\Omega})$.  For a non-negative integer $k$ and a number
$0<\alpha \leq 1$ the general H\"{o}lder space
$C^{k,\alpha} (\overline{\Omega})$ is the Banach space whose norm is defined by
    \[|y|_{k, \alpha; \Omega} = |y|_{k,0; \Omega} + [D^k y]_{ \alpha;
    \Omega}\]
where $D^k y$ denotes the $k$th order partial derivatives of $y$,
    \[|y|_{k,0; \Omega} = \sup_{\Omega} \sum_{i=0}^{k} |D^i y(x)|  \]
and
    \[ [D^k y]_{ \alpha; \Omega} =\sup_{x_1, x_2 \in \Omega}
    \frac{|D^k y(x_1)-D^k y(x_2)|}{|x_1-x_2|^{\alpha}}. \]
Notice that every  $ C^{k,\alpha} (\partial \Omega)$ function can always be
extended to a $C^{k,\alpha} (\bar{\Omega})$ function and, on
the other hand, every $C^{k,\alpha} (\bar{\Omega})$ function
can be restricted on the boundary to produce a $ C^{k,\alpha} (\partial
\Omega)$ function.  Notice also that both conversions can be carried out in a manner that preserve the norm, i.e. the corresponding $ C^{k, \alpha} (\partial \Omega)$ norm and $C^{k,\alpha} (\bar{\Omega})$ are equivalent. Based on such an observation we will not distinguish the spaces $ C^{k, \alpha}
(\partial \Omega)$ and $ C^{k, \alpha} (\bar{\Omega})$.

For the existence and uniqueness we make the following
assumptions. The definition of the allowable set of control
implies that there is a constant $M_1$ such that
    \[ |u|_{0, \alpha; \Omega}  \leq M_1.\]
We assume that $\varphi$ satisfies
\begin{enumerate}
   \item For $x\in \Omega$, $y \in \RR$
    and $|u|_{0; \Omega}  \leq M_1$, there is a positive constant $c_0$ such that
   \begin{equation} \label{phi1}
    \pd{\varphi}{y} (x,y,u) \geq c_0 >0.
    \end{equation}
    \item For $x\in \Omega$, $|y|_{0; \Omega}  \leq M_0$
    and $|u|_{0; \Omega}  \leq M_1$, there is a constant
    $M_2$ depending on $M_0$ and $M_1$ such that
         \begin{equation} \label{phi2}
          \left| \varphi (x,y,u) \right| ,  \left| D\varphi(x,y,u)
        \right|, | D^2\varphi(x,y,u)| , | D^3\varphi(x,y,u) | \leq M_2.\end{equation}
        \end{enumerate}

Let us now consider the linear problem
    \begin{equation} \label{lineqn} -\Delta y + y =0
     \quad
    \mbox{in}  \;   \Omega, \qquad
     \Delta_{\partial \Omega} y +
    \partial_{\nu} y = \alpha y +h  \quad
    \mbox{on}   \;   \partial \Omega
    \end{equation}
where $\alpha, h\in C^{\alpha} (\partial \Omega)$ are given functions such that
	\begin{equation}
 \alpha(x) \geq c_0 >0.
\label{yterm}
\end{equation}
We may assume that the $c_0$ in (\ref{yterm}) and (\ref{phi1}) are the same.
	
	 As a starting point we quote Lemma 1.1, Theorem 1.5 and Theorem 1.6 of \cite{luo:1} here.

\begin{lemma} \label{lemma1} If $x_0 \in \partial \Omega$ is a maximum point of $y \in C^{2}(\bar{\Omega})$ then at $x_0$ we have $\Delta_{\partial \Omega} y \leq 0$ and $ \partial_{\nu} y \leq 0$ so that
	\begin{equation}
	\Delta_{\partial \Omega} y +     \partial_{\nu} y \leq 0.
\end{equation}
\end{lemma}

\begin{theorem} \label{thm1} Suppose that $\Omega$ is a $C^{2, \alpha}$ domain and (\ref{yterm}) holds. If  $y \in C^{2, \alpha}$ is a
solution of (\ref{lineqn}) then
    \begin{equation} \label{apri}
    |y|_{2, \alpha; \Omega} \leq C ( |y|_{0; \Omega} + |h|_{0, \alpha; \partial
    \Omega})
    \end{equation}
where $C$ only depends on the geometry of $\partial \Omega$, $c_0$, $|\alpha|_{0, \alpha; \Omega}$ and
$n$.
\end{theorem}

\begin{theorem} \label{thm2} If $\Omega$ is a $C^{2, \alpha}$ domain then the boundary value problem (\ref{lineqn}) has a unique $C^{2, \alpha}$ solution for every $h \in C^{2, \alpha}(\partial \Omega)$.
\end{theorem}

To handle nonlinear problems we frequently need the following fact:
\begin{lemma}\label{compose}
Suppose $\psi : \Omega \times \RR^k \rightarrow \RR$ is a $C^1$ function satisfying
	\[ |\psi|, |D\psi | \leq M_2 \]
for a constant $M_2$.  If $u_i(x) \in  C^{\alpha}( \bar{\Omega})$ for $i=1,2,\ldots, k$ then
$\psi(x, u_1(x), \ldots, u_k(x)) \in  C^{\alpha}( \bar{\Omega})$ and
\begin{equation}
 |\psi|_{0, \alpha; \Omega} \leq M_2 \left(1+ d^{1-\alpha} + \sum_{i=1}^{k} |u_i|_{0, \alpha; \Omega})\right)
\label{eqcompose}
\end{equation}
where $d$ is the diameter of $\Omega$.
\end{lemma}

\pro
Obviously
	\[ |\psi|_{0;\Omega} \leq M_2.\]
Also by the assumption on $\psi$ we have
	 \begin{eqnarray*} [\psi]_{\alpha ;\Omega} & = & \sup_{x_1, x_2 \in \Omega}
    \frac{|\psi ( x_1, u_1(x_1),\ldots,
    u_k(x_1))-\psi ( x_2, u_1(x_2),\ldots,
    u_k(x_2))|}{|x_1-x_2|^{\alpha}} \\
    & \leq & \sup_{x_1, x_2 \in \Omega}
    \left( \left| \pd{\psi}{x} \right||x_1-x_2|^{1-\alpha} +
    \sum_{i=1}^{k}  \left| \pd{\psi}{u_i} \right|  \frac{|u_i( x_1)-u_i(x_2)|}{|x_1-x_2|^{\alpha}}
    \right)  \\
    & \leq & M_2 \left( d^{1-\alpha} + \sum_{i=1}^{k}[u_i]_{\alpha ;\Omega}\right).
    \end{eqnarray*}
Adding these up yields (\ref{eqcompose}).

Now we turn to our nonlinear problem (\ref{state}).  Suppose that for each $u \in U$ the problem (\ref{state}) has a classical solution $y$. We show first that the $C^{0}( \Omega)$ norm of $y$ is bounded by the $C^{0}( \Omega)$ norm of $u$ and therefore it is bounded as $U$ is a bounded set in $C^{0}( \Omega)$.
Let $M_0 = |y|_{0; \Omega}
=\sup_{\Omega} |y | \geq 0$. Without loss of generality we assume
$ M_0 = y(x_0)$ for some $x_0 \in \bar{\Omega}$.  By the {\em weak maximum principle} we have $x_0
\in
\partial \Omega$. Then Lemma \ref{lemma1} implies
    \begin{eqnarray*}
     0 & \geq & \Delta_{\partial \Omega} y ( x_0)+
    \partial_{\nu} y( x_0) = \varphi ( x_0, M_0, u(x_0)) \\
    & = & M_0 \int_{0}^{1} \pd{\varphi}{y}
    ( x_0, tM_0, u(x_0))\; dt + \varphi ( x_0, 0, u(x_0)) \\
    & \geq & c_0 M_0 + \varphi ( x_0, 0, u(x_0)).
    \end{eqnarray*}
Thus by (\ref{phi1}) and (\ref{phi2}) we have the desired bound for $y$:
    \begin{equation} \label{ybound}
    M_0 \leq \frac{1}{c_0} \sup_{\Omega} |\varphi ( \cdot, 0, u(\cdot)) |
    \leq \frac{1}{c_0} M_2.
    \end{equation}
Next we show that the $C^{2,\alpha}( \Omega)$ norm of $y$ is also bounded.  For this purpose we put $h(x) = \varphi ( x, y(x), u(x))$.  For convenience we may assume, from now on, that $d\geq 1$ so $d^{1-\alpha} \leq d$.  Then, by Lemma \ref{compose}, we have
      \begin{equation}\label{hbound}
    |h|_{0, \alpha; \partial
    \Omega} \leq M_2 ( 1+d + M_1+|y|_{0,\alpha ;\Omega} ).
    \end{equation}
By inserting (\ref{ybound}) and (\ref{hbound}) into (\ref{apri})
we obtain
    \begin{equation}\label{apri1}
     |y|_{2, \alpha; \Omega} \leq C ( M_0 + M_2 ( 1+d + M_1+|y|_{0,\alpha ;\Omega}
     )) = C_1 + C_2 |y|_{0,\alpha ;\Omega}
    \end{equation}
for some constants $C_1$ and $C_2$.
Now we recall the well known interpolation inequality (Lemma 6.35
of \cite{GilbTrud:1}):
	\begin{equation}
	|y|_{0,\alpha ;\Omega} \leq C |y|_{0 ;\Omega} +\varepsilon |y|_{2, \alpha ;\Omega}
	\leq  \frac{1}{c_0} C_3M_2 + \varepsilon |y|_{2, \alpha ;\Omega}
\label{interpo}
\end{equation}
for some constant $C=C(\varepsilon, \Omega)$.  By choosing $\varepsilon = 1/(2C_2)$ we have proved
the important \emph{a priori} estimate for the solutions of (\ref{state}):

\begin{theorem} \label{bound} Assume that $\varphi $ satisfies (\ref{phi1}) and (\ref{phi2}), and $u \in U$.
If $y$ is a solution of (\ref{state}) then we have
    \[ |y|_{2, \alpha; \Omega} \leq C\]
for a constant $C$ depending on $M_1$, $M_2$, $\Omega$ and $n$.
\end{theorem}

With $u$ fixed, for all $z \in   C^{2,\alpha} (\partial \Omega)$
we define the operator $T$ by letting $y=Tz$ be the unique
solution in $ C^{2,\alpha} (\partial \Omega)$ of the linear problem
\begin{eqnarray} \label{linear} -\Delta y + y =0
     \qquad &
    \mbox{in}  & \quad  \Omega, \nonumber \\
     \Delta_{\partial \Omega} y +
    \partial_{\nu} y = \varphi (\cdot, z,u)  \quad &
    \mbox{on}  & \quad   \partial \Omega.
    \end{eqnarray}
The unique solvability of this linear problem is guaranteed by Theorem \ref{thm2}.  By applying the Leray-Schauder fixed point theorem together with Theorem \ref{bound} to the operator $T$ we
obtain the existence of a fixed point of $T$ which is then a
solution of our state equation (\ref{state}).

To see whether the solution of (\ref{state}) is unique we suppose
that $y_1$ and $y_2$ are two solutions.  Then the difference
$y=y_1-y_2$ will satisfy
 \begin{eqnarray} \label{uniq} -\Delta y + y =0
     \qquad &
    \mbox{in}  & \quad  \Omega, \nonumber \\
     \Delta_{\partial \Omega} y +
    \partial_{\nu} y = \varphi (\cdot, y_1,u) -\varphi (\cdot, y_2,u) \quad &
    \mbox{on}  & \quad   \partial \Omega.
    \end{eqnarray}
If $y_1 \neq y_2$ we may assume $M=y_1(x_0)-y_2(x_0) =
\sup_{\Omega} |y_1-y_2| >0$.  The weak maximum principle
implies that $x_0 \in \partial \Omega$.  Then by Lemma \ref{lemma1}, at $x_0$ we have
    \[ 0 \geq \Delta_{\partial \Omega} y +
    \partial_{\nu} y = \varphi (\cdot, y_1,u) -\varphi (\cdot,
    y_2,u) >0,\]
a contradiction.  In summary we have
\begin{theorem} \label{unique} Assume that $\varphi $ satisfies (\ref{phi1}) and (\ref{phi2}), and $u \in U$. Then the state equation (\ref{state}) has a unique
solution $y  \in C^{2,\alpha} (\bar{\Omega})$.
\end{theorem}

\section{Differentiability}
\label{sect:3}
Given $u\in C^{\alpha} (\partial \Omega)$ we denote by $y_u$ the
solution of the boundary value problem (\ref{state}). This
correspondence defines a mapping $G: C^{\alpha} (\partial \Omega)
\rightarrow C^{2, \alpha} (\bar{\Omega})$ by $y_u=G(u)$.

\begin{theorem} The mapping $y=G(u)$ is Fr\'{e}chet
differentiable.  Let $G'(u)$ 
be the Fr\'{e}chet
derivative of $G$ and $z=\langle G'(u), v\rangle$ where $v\in
C^{\alpha} (\partial \Omega)$.  Then $z$ is the unique solution of the
boundary value problem
    \begin{eqnarray} \label{gdash1} -\Delta z + z =0
     \qquad &
    \mbox{in}  & \quad  \Omega,  \\
    \label{gdash2} \Delta_{\partial \Omega} z +
    \partial_{\nu} z = \frac{\partial \varphi }{\partial y}
    (\cdot, y, u)  z + \frac{\partial \varphi }{\partial u}
    (\cdot, y, u) v  \quad &
    \mbox{on}  & \quad   \partial \Omega.
    \end{eqnarray}
\end{theorem}

\pro  We first prove that $G$ is Gateaux-differentiable and calculate the
G-derivative  $dG(u)$. Let $v\in C^{\alpha} (\partial \Omega)$ and
consider $y_t=G(u+tv)$ and $y=G(u)$.  It follows that
     \begin{eqnarray} \label{eqn1} -\Delta y_t + y_t =0
     \qquad &
    \mbox{in}  & \quad  \Omega, \nonumber \\
     \Delta_{\partial \Omega} y_t +
    \partial_{\nu} y_t = \varphi (\cdot, y_t,u+tv )  \quad &
    \mbox{on}  & \quad   \partial \Omega
    \end{eqnarray}
and
     \begin{eqnarray} \label{eqn2} -\Delta y + y =0
     \qquad &
    \mbox{in}  & \quad  \Omega, \nonumber \\
     \Delta_{\partial \Omega} y +
    \partial_{\nu} y = \varphi (\cdot, y,u)  \quad &
    \mbox{on}  & \quad   \partial \Omega.
    \end{eqnarray}
By subtracting (\ref{eqn2}) from (\ref{eqn1}) we see that $w_t=
y_t-y$ satisfies
    \begin{eqnarray} \label{eqnwt} -\Delta w_t + w_t =0
     \qquad &
    \mbox{in}  & \quad  \Omega, \nonumber \\
     \Delta_{\partial \Omega} w_t +
    \partial_{\nu} w_t = \varphi (\cdot, y_t,u+tv )- \varphi (\cdot, y,u)  \quad &
    \mbox{on}  & \quad   \partial \Omega.
    \end{eqnarray}
We can assume that $t$ is bounded, say $|t| \leq 1$.  Theorem \ref{bound}  guarantees that the $C^{2,\alpha}( \Omega)$ norms of both $y$ and $y_t$ are uniformly bounded. Once $u$ and $v$ are chosen their norms are also independent of $t$.  We now want to show that the norm $| w_t |_{2,\alpha;\Omega} $ is also uniformly bounded with respect to $|t|<1$. For this purpose we consider the composite function $h(x)=\varphi (x, y_t(x),u(x)+tv(x) )-\varphi  (x, y(x), u(x))$.  By Lemma \ref{compose} again we have
		\begin{equation}
 	| h |_{0,\alpha;\Omega} := | h |_{0;\Omega} + [ h ]_{\alpha;\Omega} \leq C
\label{hbound1}
\end{equation}
for a constant $C$ that is independent of $t$.  Using Theorem \ref{thm1} we fist obtain
	\begin{equation}
| w_t |_{2,\alpha;\Omega} \leq C_1( |w_t|_{0; \Omega} + |h|_{0, \alpha; \partial
    \Omega}) \leq  C_1( |w_t|_{0; \Omega} + C)
\label{wt2a}
\end{equation}
Then the problem is reduced to the estimation of $|w_t|_{0; \Omega} $.  The previous argument does not work here because equation (\ref{eqnwt}) does not have the $w_t$ term.  However this situation is covered by Lemma 1.4 of \cite{luo:1} which gives
	\begin{equation}
|w_t|_{0; \Omega} \leq C_2 \sup_{x \in \Omega}
	 |\varphi (x, y_t(x),u(x)+tv(x) )-\varphi  (x, y(x), u(x)) | \leq 2 C_2 M_2
\label{wt0}
\end{equation}
where the constant $C_2$ is independent of $t$.  A substitution of (\ref{wt0}) into (\ref{wt2a}) produces
\begin{equation}
| w_t |_{2,\alpha;\Omega} \leq C_3
\label{wtunif}
\end{equation}
for a constant $C_3$ independent of $t$, that is, $| w_t |_{2,\alpha;\Omega} $ is uniformly bounded.
Therefore, up to a subsequence, the following limits exist in
$C^{2,\alpha} (\bar{\Omega})$:
    \[  \lim_{t \rightarrow 0} w_t = w \qand  \lim_{t \rightarrow 0} y_t
    = \lim_{t \rightarrow 0} (w_t+y) = w +y.\]
By taking limit in (\ref{eqnwt}) as $t \rightarrow 0$ we have
\begin{eqnarray} \label{eqnw} -\Delta w + w =0
     \qquad &
    \mbox{in}  & \quad  \Omega, \nonumber \\
     \Delta_{\partial \Omega} w +
    \partial_{\nu} w = \varphi (\cdot, y+w,u)- \varphi (\cdot, y,u)  \quad &
    \mbox{on}  & \quad   \partial \Omega.
    \end{eqnarray}
Since $\varphi$ is increasing in $y$ variable, the only solution
satisfying (\ref{eqnw}) is $w=0$.  From this we conclude $ \lim_{t
\rightarrow 0} y_t = y$.

Now consider $z_t=w_t/t$. Dividing  (\ref{eqnwt}) by $t$ yields
    \begin{eqnarray} \label{eqnzt} -\Delta z_t + z_t =0
     \qquad &
    \mbox{in}  & \quad  \Omega, \nonumber \\
     \Delta_{\partial \Omega} z_t +
    \partial_{\nu} z_t = \alpha_t z_t + \beta_t v  \quad &
    \mbox{on}  & \quad   \partial \Omega
    \end{eqnarray}
where
    \[ \alpha_t = \int_{0}^{1} \frac{\partial \varphi }{\partial y}
    (x, \tau y_t+(1-\tau)y, u+tv)\; d\tau \qand \beta_t =\int_{0}^{1}
    \frac{\partial \varphi }{\partial u}
    (x, y, u+\tau tv)\; d\tau.\]
Obviously  $\alpha_t \in C^{2,\alpha} (\bar{\Omega})$.  It follows from Lemma \ref{compose} that
	\[ | \alpha |_{0,\alpha;\Omega}  \leq  M_2(d+|y |_{0,\alpha;\Omega}+ | y_t |_{0,\alpha;\Omega}+|u |_{0,\alpha;\Omega}+|v |_{0,\alpha;\Omega}) \leq C_4
		\]
for a constant $C_4$. Notice also that $\alpha_t \geq c_0$ and $|\beta_t |\leq M_2$ and hence Theorem \ref{thm1} implies
	\[ | z_t |_{2,\alpha;\Omega} \leq C_5(|z_t|_{0; \Omega} +|\beta_t v |_{0,\alpha;\Omega})
	\leq C_5(|z_t|_{0; \Omega} +C_6)\]
for some constants $C_5$ and $C_6$.  Finally the uniform bound for $| z_t |_{2,\alpha;\Omega} $ comes from the estimate
	\[ |z_t|_{0; \Omega} \leq \frac{1}{c_0} M_2 |v|_{0; \Omega}, \]
which is a consequence of Lemma \ref{lemma1}.
In summary we have
\begin{equation}
| z_t |_{2,\alpha;\Omega} \leq C_7
\label{ztunif}
\end{equation}
for a constant $C_7$ that is independent of $t$.
This implies that, up to a subsequence, $z_t$
converges to a function $z$ in $C^{2,\alpha} (\bar{\Omega})$
as $t \rightarrow 0$ and
    \[ \lim_{t \rightarrow 0} \alpha_t = \frac{\partial \varphi }{\partial y}
    (x, y, u) \qand \lim_{t \rightarrow 0} \beta_t = \frac{\partial \varphi }{\partial u}
    (x, y, u).\]
Taking limit in (\ref{eqnzt}) gives
    \begin{eqnarray} \label{eqnz} -\Delta z + z =0
     \qquad &
    \mbox{in}  & \quad  \Omega, \nonumber \\
     \Delta_{\partial \Omega} z +
    \partial_{\nu} z = \frac{\partial \varphi }{\partial y}
    (\cdot, y, u)  z + \frac{\partial \varphi }{\partial u}
    (\cdot, y, u) v  \quad &
    \mbox{on}  & \quad   \partial \Omega
    \end{eqnarray}
which means that $z=\langle dG(u), v\rangle$ is the solution of
(\ref{gdash1}) and (\ref{gdash2}).

The uniqueness of $z$ is guaranteed by Theorem \ref{thm2} as (\ref{eqnz}) is a linear equation.

Next we examine the continuity of $dG$.  Notice that $dG(u)\in
{\cal L} ( C^{\alpha} (\partial \Omega),  C^{2, \alpha}
(\bar{\Omega}))$ and
    \[ \| dG(u) \| = \sup_{\| v \| =1} | \langle dG(u), v\rangle |_{2,\alpha;\Omega}. \]
Therefore to prove the continuity of $dG(u)$ is to prove that as $u_1 \rightarrow u$ in $C^{\alpha} (\bar{\Omega})$
\[ \| dG(u_1)- dG(u)\| = \sup_{\| v \| =1} | \langle dG(u_1), v\rangle -
\langle dG(u), v\rangle |_{2,\alpha;\Omega} \rightarrow 0. \]
For any $v \in C^{\alpha} (\partial \Omega)$ with $\| v \| =|v|_{\alpha;\Omega}=1$
consider
    $ z_1= \langle dG(u_1), v\rangle$ and $z= \langle dG(u), v\rangle$.  Then we know that $w_1=z_1-z$
is a solution of
    \begin{eqnarray} \label{contin} -\Delta w_1 + w_1=0
     \qquad &
    \mbox{in}  & \quad  \Omega, \nonumber \\
     \Delta_{\partial \Omega} w_1 +
    \partial_{\nu} w_1 = \frac{\partial \varphi }{\partial y}
    (\cdot, y, u_1)  z_1 + \frac{\partial \varphi }{\partial u}
    (\cdot, y, u_1) v & & \nonumber \\
     - \frac{\partial \varphi }{\partial y}
    (\cdot, y, u)  z - \frac{\partial \varphi }{\partial u}
    (\cdot, y, u) v \quad &
    \mbox{on}  & \quad   \partial \Omega.
    \end{eqnarray}
All we need to show is that $w_1 \rightarrow 0$ in $C^{2,\alpha} (\bar{\Omega})$ uniformly with respect to $|v|_{\alpha;\Omega}=1$, as $u_1 \rightarrow u$ in $C^{\alpha} (\bar{\Omega})$.  To this end we rewrite equation (\ref{contin}) in the form
	\begin{eqnarray} \label{contin1} -\Delta w_1 + w_1=0
     \qquad &
    \mbox{in}  & \quad  \Omega, \nonumber \\
     \Delta_{\partial \Omega} w_1 +
    \partial_{\nu} w_1 = \sigma w_1 + \gamma \quad &
    \mbox{on}  & \quad   \partial \Omega
    \end{eqnarray}
where
\[ \sigma = \frac{\partial \varphi }{ \partial y}
    (x, y, u_1) \]
 and
 \[ \gamma =\left( \frac{\partial \varphi }{\partial y}
    (\cdot, y, u_1)- \frac{\partial \varphi }{\partial y}
    (\cdot, y, u)\right)  z + \left(\frac{\partial \varphi }{\partial u}
    (\cdot, y, u_1) - \frac{\partial \varphi }{\partial u}
    (\cdot, y, u) \right)v.\]
If we put
	\[ A =  \int_{0}^{1} \frac{\partial^2 \varphi }{\partial u \partial y}
    (x, y, u+t(u_1-u))\; dt \qand B =  \int_{0}^{1} \frac{\partial^2 \varphi }{\partial u^2}
    (x, y, u+t(u_1-u))\; dt\]
then $\gamma$ can be written as
\[ \gamma =(Az+Bv)(u-u_1).\]
From the assumption on $\varphi$ we know that $(Az+Bv) \in  C^{\alpha} (\bar{\Omega})$ and hence
	\[ |\gamma|_{\alpha;\Omega} \leq C_8 |u_1-u|_{\alpha;\Omega}.\]
By Theorem \ref{bound} we then have
    \[ |w_1|_{2,\alpha;\Omega}
    \leq C_9( |w_1|_{0;\Omega} +|\gamma|_{\alpha;\Omega}) \leq C_9( |w_1|_{0;\Omega} +C_8 |w_1|_{\alpha;\Omega}) \leq C_9(1+C_8)|w_1|_{\alpha;\Omega}.
    \]
Now the continuity of  $dG(u)$ follows because $|w_1|_{\alpha;\Omega}\rightarrow 0 $ uniformly with respect to $|v|_{\alpha;\Omega}=1$. Finally, since $G(u)$ is continuously Gateaux differentiable,  we conclude that $G(u)$ is also Fr\'{e}chet differentiable and that the Fr\'{e}chet derivative $G'(u)$ is equal to $dG(u)$.

Now we are in the position to establish the differentiability of the objective functional
$J(u)$.

\begin{theorem} \label{diffj} The functional $J$ is Fr\'{e}chet differentiable and
for every $u, v \in   C^{\alpha} (\partial \Omega)$ and $y = G(u)$
we have
    \[ \langle J'(u), v\rangle =\int_{\partial \Omega} \left[ \pd{g}{u}(s, y,u)- \frac{\partial \varphi }{\partial u}
    (s, y, u)w \right] v \; ds \]
where $w$ is the solution of
    \begin{eqnarray} \label{w1} -\Delta w + w =
    \pd{f}{y}(\cdot, y)
     \qquad &
    \mbox{in}  & \quad  \Omega,  \\
     \label{w2}  \Delta_{\partial \Omega} w +
    \partial_{\nu} w = \frac{\partial \varphi }{\partial y}
    (\cdot, y, u)  w - \pd{g}{y}(\cdot,y, u) \quad &
    \mbox{on}  & \quad   \partial \Omega.
    \end{eqnarray}
\end{theorem}

\pro Define
    \[ H(y,u) =\int_{\Omega} f(x, y(x)) \; dx
    +\int_{\partial \Omega} g(s, y(s),u(s)) \; ds.
     \]
It follows that
    \[ J(u) = H(G(u),u).\]
It is obvious that $H$ is differentiable and for every $\bar{y}$
and $\bar{u}$ we have
    \[ \langle \pd{H}{y}(y,u), \bar{y}\rangle = \int_{\Omega} \pd{f}{y}(x, y(x)) \bar{y} \; dx
    +\int_{\partial \Omega} \pd{g}{y}(s, y(s),u(s)) \bar{y}\; ds
     \]
and
  \[ \langle \pd{H}{u}(y,u), \bar{u}\rangle = \int_{\partial \Omega} \pd{g}{u}(s, y(s),u(s)) \bar{u}\;
  ds.
     \]
By the chain rule we have
    \[ \langle J'(u), v \rangle = \langle \pd{H}{y}(y,u)G'(u)+  \pd{H}{u}(y,u), v\rangle \]
which then gives
    \begin{eqnarray}\label{jdash1}
     \langle J'(u), v\rangle &=&  \int_{\Omega} \pd{f}{y}(x, y(x))G'(u)v \; dx
      + \int_{\partial \Omega} \pd{g}{y}(s, y(s),u(s))G'(u)v \; ds
      \nonumber \\
    & &    + \int_{\partial \Omega} \pd{g}{u}(s, y(s),u(s))v \; ds  \nonumber \\
    & = &  \int_{\Omega} \pd{f}{y}(x, y(x))z(x) \; dx
      + \int_{\partial \Omega} \pd{g}{y}(s, y(s),u(s))z(s) \; ds
      \nonumber \\
    & &    + \int_{\partial \Omega} \pd{g}{u}(s, y(s),u(s))v \; ds
    \end{eqnarray}
where $z$ is the solution of (\ref{gdash1}) and (\ref{gdash2}).

Let $w$ be the solution of (\ref{w1}) and (\ref{w2}).  Subtracting
$w$ times (\ref{gdash1}) from $z$ times (\ref{w1}) and applying
the Green's second identity yields
    \begin{eqnarray*}  \int_{\Omega} \pd{f}{y}(x, y)z \; dx &=&
     \int_{\Omega}( w\Delta z-z \Delta w)  dx =
     \int_{\partial \Omega} (z \partial_{\nu} w-  w\partial_{\nu} z)
     \;  ds \\
     & = & \int_{\partial \Omega} \left\{ z \partial_{\nu} w -w[-\Delta_{\partial \Omega} z +
     \frac{\partial \varphi }{\partial y}
    (s, y, u)  z + \frac{\partial \varphi }{\partial u}
    (s, y, u) v] \right\} \;  ds \\
     & = & \int_{\partial \Omega}  w\Delta_{\partial \Omega} z \; ds \\
     & & +
     \int_{\partial \Omega} \left\{ [ \partial_{\nu} w- \frac{\partial \varphi }{\partial y}
    (s, y, u) w   ]z - \frac{\partial \varphi }{\partial u}
    (s, y, u) vw \right\}  ds.
     \end{eqnarray*}
On the manifold $\partial \Omega$ there holds a boundary version of Green's identity, see Lemma 16.1 of \cite{GilbTrud:1}.  Using this boundary version of Green's identity we have
    \[ \int_{\partial \Omega}  w\Delta_{\partial \Omega} z \; ds =
    \int_{\partial \Omega}  z\Delta_{\partial \Omega} w \; ds\]
and hence,
    \begin{eqnarray}  \int_{\Omega} \pd{f}{y}(x, y)z \; dx &=&
          \int_{\partial \Omega} \left\{ [\Delta_{\partial \Omega} w+
           \partial_{\nu} w- \frac{\partial \varphi }{\partial y}
    (s, y, u) w   ]z - \frac{\partial \varphi }{\partial u}
    (s, y, u) wv \right\}  ds \nonumber \\
    & = & \int_{\partial \Omega} \left\{ -\pd{g}{y}(s, y,u) z - \frac{\partial \varphi }{\partial u}
    (s, y, u)w v \right\}  ds .
     \end{eqnarray}
A substitution of this into (\ref{jdash1}) gives
    \[ \langle J'(u), v\rangle  = \int_{\partial \Omega} \left[ \pd{g}{u}(s, y,u)- \frac{\partial \varphi }{\partial u}
    (s, y, u)w\right] v \; ds.
    \]

\section{Main result}
\label{sect:4}
A function $\bar{u } \in U$  is said to be a local solution, or a
locally optimal control, of (\ref{prob}) if there is a number
$\delta > 0$ such that $J(u) \geq J(\bar{u})$ holds for all $u \in
U$ satisfying $|u-\bar{u}| < \delta$, with their associated state
$y $ and the state constraint on $y$.  Our main result is the first order necessary condition for a $\bar{u } \in U$ to be a local solution.

\begin{theorem} If $\bar{u } \in U$  is a local solution of (\ref{prob})
then there exist a real number $\lambda \geq 0$, a function
$\bar{y} \in C^{2,\alpha} (\bar{\Omega})$ and a function
$\bar{w} \in C^{2,\alpha} (\bar{\Omega})$ such that $\bar{y}$
satisfies
   \begin{eqnarray} \label{nec1} -\Delta \bar{y} + \bar{y} =
    0
     \qquad &
    \mbox{in}  & \quad  \Omega,  \nonumber \\
       \Delta_{\partial \Omega} \bar{y} +
    \partial_{\nu} \bar{y} =  \varphi (\cdot, \bar{y}, \bar{u}) \quad &
    \mbox{on}  & \quad   \partial \Omega
    \end{eqnarray}
$\bar{w}$ satisfies
     \begin{eqnarray} \label{nec2} -\Delta \bar{w} + \bar{w} =
    \pd{f}{y}(\cdot, \bar{y}) + \lambda \pd{a}{y}(\cdot, \bar{y})
     \qquad &
    \mbox{in}  & \quad  \Omega,  \nonumber \\
      \Delta_{\partial \Omega} \bar{w} +
    \partial_{\nu} \bar{w} = \frac{\partial \varphi }{\partial y}
    (\cdot, \bar{y}, \bar{u})  \bar{w} - \pd{g}{y}(\cdot, \bar{y}, \bar{u}) -\lambda
    \pd{b}{y}(\cdot, \bar{y}, \bar{u}) \quad &
    \mbox{on}  & \quad   \partial \Omega
    \end{eqnarray}
and
     \begin{equation} \label{nec3}
     \int_{\partial \Omega} \left[ \pd{g}{u}(s, \bar{y},\bar{u})-
     \frac{\partial \varphi }{\partial u}
    (s, \bar{y}, \bar{u})\bar{w}\right](u- \bar{u}) \; ds \geq 0
    \end{equation}
for all  $u \in U$.
\end{theorem}

\pro  Let $y = G(u)$ be the solution of
(\ref{nec1}) corresponding to $u$ and $\rho (u) = F(G(u))$ where
$F$ is the constraint functional in our optimal control problem
(\ref{prob}). By the theorem of Lagrange multiplier there is a
$\lambda \geq 0$ such that
     \begin{equation} \label{lag}
     \langle J'(\bar{u})+\lambda \rho '(\bar{u}), u- \bar{u} \rangle
     \geq 0
    \end{equation}
for all  $u \in U$.  In order to obtain the necessary conditions
stated in the theorem the only thing remaining is to calculate the
derivative $\rho '(\bar{u})$ of the constrain functional $F$. As
in the proof of Theorem \ref{diffj} we have
    \begin{eqnarray}\label{psidash1}
     \langle \rho '(\bar{u}), v\rangle &=&  \int_{\Omega} \pd{a}{y}(x, \bar{y}(x))G'(\bar{u})v \; dx
      + \int_{\partial \Omega} \pd{b}{y}(s, \bar{y}(s))G'(\bar{u})v \; ds
      \nonumber \\
       & = &  \int_{\Omega} \pd{a}{y}(x, \bar{y}(x))z(x) \; dx
      + \int_{\partial \Omega} \pd{b}{y}(s, \bar{y}(s))z(s) \; ds
        \end{eqnarray}
where $z$ is the solution of (\ref{gdash1}) and (\ref{gdash2})
corresponding to $\bar{y}$ and $\bar{u}$.

Let $w_1$ be the solution of
\begin{eqnarray} \label{constrain1} -\Delta w_1 + w_1 =
    \pd{a}{y}(\cdot, \bar{y})
     \qquad &
    \mbox{in}  & \quad  \Omega,  \nonumber \\
      \Delta_{\partial \Omega} w_1 +
    \partial_{\nu} w_1 = \frac{\partial \varphi }{\partial y}
    (\cdot, \bar{y}, \bar{u})  w_1 - \pd{b}{y}(\cdot, \bar{y}, \bar{u}) \quad &
    \mbox{on}  & \quad   \partial \Omega.
    \end{eqnarray}
By the same argument as in the proof of Theorem \ref{diffj} we
obtain
    \[ \langle \rho '(\bar{u}), v\rangle  = \int_{\partial \Omega} \left[ - \frac{\partial \varphi }{\partial u}
    (s, \bar{y}, \bar{u})w_1\right] v \; ds.
    \]
Finally we put $\bar{w}= w + \lambda w_1$ where $w $ is the
solution of (\ref{w1}) and (\ref{w2}) corresponding to $\bar{y}$
and $\bar{u}$.  Then $\bar{w}$ is the solution of (\ref{nec2}) and
    \[  \langle J'(\bar{u})+\lambda \rho '(\bar{u}), u- \bar{u} \rangle
    =  \int_{\partial \Omega} \left[ \pd{g}{u}(s, \bar{y},\bar{u})-
     \frac{\partial \varphi }{\partial u}
    (s, \bar{y}, \bar{u})\bar{w}\right](u- \bar{u}) \; ds.\]

\section{Remarks}
\label{sect:5}
A second order necessary optimality condition for the problem (\ref{prob}) can be easily
established using exactly the same arguments as in Section \ref{sect:3} and \ref{sect:4}.
To keep this article short we leave the derivation and formulation for such a result to the interested readers.

As mentioned in the Introduction the initial-boundary value problem of a parabolic equation with a parabolic Venttsel boundary condition arises in the engineering problem of heat conduction.  When one considers the optimal control problem (\ref{prob}) with such a state equation, similar optimality conditions are expected.   The study on such a problem is currently undergoing and the result will be published in the near future.  Here we just point out the formulation of the state equation under the consideration:
\begin{eqnarray} \label{statep} \pd{y}{t}-\Delta y + y =0
     \qquad &
    \mbox{in}  & \quad  \Omega, \nonumber \\
     \pd{y}{t} -\Delta_{\partial \Omega} y -
    \partial_{\nu} y = \varphi (\cdot, y,u)  \quad &
    \mbox{on}  & \quad   \partial \Omega.
    \end{eqnarray}
The theoretical frame work in dealing with such a problem has been well established. This includes existence and uniqueness of the solution in a proper function space, as well as the \emph{a priori} estimates.  For details, see \cite{apush:2} and \cite{luo:3} and the references therein.

Finally we would like to point out that there is no difficulty to extend the results in this paper to the case when the state equation is a general second order elliptic equation with a general Venttsel boundary condition:
\begin{eqnarray} \label{stateg} a^{ij} D_{ij}y + b^i D_i y + c y =0
     \qquad &
    \mbox{in}  & \quad  \Omega, \nonumber \\
     \alpha^{ij}\partial_{ij} y+ \partial_{\nu} y = \varphi (\cdot, y,u)  \quad &
    \mbox{on}  & \quad   \partial \Omega
    \end{eqnarray}
where  $a^{ij}, b^i, c, \alpha^{ij}$ are all constants
satisfying the conditions:
\begin{description}
\item[(i)]  $\{a^{ij}\}$ is a positive definite symmetric matrix
with the smallest eigenvalue $\sigma
>0$ and $c < 0$;
\item[(ii)]  $\{\alpha^{ij}\}$ is a positive definite symmetric matrix
with the smallest eigenvalue $\gamma
>0$.
\end{description}
In this general case, when proving a similar result to Theorem \ref{diffj}, due to the lack of Green's second identity, the result will take a more complicated form.




\end{document}